\documentstyle[12pt]{article} 
\newcommand{\doublespace}{
   \renewcommand{\baselinestretch}{1.2}
   \large\normalsize}

\def \Z{\Bbb Z}
\def \C{\Bbb C}

\def \mod{{\rm mod\,}}
\def \<{\langle} 
\def \>{\rangle} 
\def \t{\tau }

\def \e{\epsilon }
\def \l{\lambda }

\def \g{\gamma}

\def \qed{\mbox{ $\square$}}
\def \pf {\noindent {\bf Proof:} \,}
\input amssym.def
\input amssym
\doublespace
\begin{document}
\newtheorem{thm}{Theorem}[section]
\newtheorem{th}{Theorem}
\newtheorem{prop}[thm]{Proposition}
\newtheorem{coro}[thm]{Corollary}
\newtheorem{lem}[thm]{Lemma}
\newtheorem{rem}[thm]{Remark}
\newtheorem{de}[thm]{Definition}
\begin{center}
{\Large {\bf Transformation Laws for Theta functions}} \\
\vspace{0.5cm}
Chongying Dong\footnote{Supported by NSF grant 
DMS-9700923 and a research grant from the Committee on Research, UC Santa Cruz.} and Geoffrey Mason\footnote{Supported by NSF grant DMS-9700909 and a research grant from the Committee on Research, UC Santa Cruz.}\\
Department of Mathematics, University of
California, Santa Cruz, CA 95064
\end{center}

\section{Introduction} 

We prove some results which extend the classical 
theory, due to Hecke-Schoeneberg [H], [S1], of the transformation laws
of theta functions. Although our results are classical in natural, they were
suggested by recent work involving modular-invariance in conformal 
field theory [DMN], and we shall say more about these connection
in the last section of the present paper. 

Let $Q$ be a positive-definite, integral quadratic form of even rank 
$f=2r$ with theta-function
\begin{equation}\label{1.1}
\theta(Q,\tau)=\sum_{m\in\Z^f}e^{2\pi iQ(m)\tau}.
\end{equation}
One knows [S2] that $\theta(Q,\tau)$ is a modular form of weight $r$ and character $\e$ on the group $\Gamma_0(N).$ Here $N$ is the level of $Q$ and $\e$  is
the Dirichlet character given by the Jacobi symbol 
$$\e(n)=\left( \frac{(-1)^{r}\det A}{n}\right)$$ for $n>0,$ with $A$ a
Gram matrix of the bilinear form $\<,\>$ corresponding to $Q.$

We fix a vector $v\in \C^f$ and define 
\begin{equation}\label{1.2}
\theta(Q,v,k,\tau)=\sum_{m\in\Z^l}\<v,m\>^ke^{2\pi i Q(m)\t}
\end{equation}
where $k$ is a nonnegative integer. Obviously, $\theta(Q,v,k,\tau)$ 
identically zero if $k$ is odd, and coincides with $\theta(Q,\t)$ if
$k=0.$ 

We define 
 \begin{equation}\label{1.3}
\Theta(Q,v,\tau,X)=\sum_{n\geq 0}\frac{2^n\theta(Q,v,2n,\tau)}{(2n)!}(2\pi i X)^n
\end{equation}
regarding this as a function on ${\frak h}\times \C$ ($\frak h$ is
the complex upper half-plane) for fixed $Q$ and $v.$ 

Our main result may then be stated as follows:
\begin{th}\label{t1} $\Theta(Q,v,\t,X)$ satisfies the following
transformation law for $\gamma=\left(
\begin{array}{cc} a & b\\ c & d\end{array}\right)$ in $\Gamma_0(N):$

\begin{equation}\label{1.4}
\Theta(Q,v,\gamma\tau,\frac{X}{(c\t+d)^2})=\e(d)(c\tau+d)^{r}\exp\left(\frac{c\<v,v\>X}{c\tau+d}\right)\Theta(Q,v,\tau,X).
\end{equation}
\end{th}

It is evident from (1.2) that scaling $v$ (i.e., replacing $v$ by
$\l v$ for a scalar $\l$) simply multiplies $\theta(Q,v,k,\t)$ by $\l^k.$ For
this reason, there are essentially only two cases, namely

(a) $v$ is a null vector i.e., $\<v,v\>=0,$ or

(b) $v$ is a unit vector i.e., $\<v,v\>=1.$

Suppose first that $v$ is a null vector. Then Hecke [H] proved that
the function $P_k(x)=\<v,x\>^k$ is a spherical harmonic of degree $k$
with respect to $Q,$ moreover every spherical harmonic of degree $k$
is a linear combination of such functions. Thus in this case (\ref{1.2})
is simply a theta function with spherical harmonic $\theta(Q,P_k,\tau),$
and 
\begin{equation}\label{1.5}
\Theta(Q,v,\tau,X)=\sum_{n\geq 0}\frac{2^n\theta(Q,P_{2n},\tau)}{(2n)!}(2\pi i X)^n.
\end{equation}
The transformation law (\ref{1.4}) then says exactly that
$\theta(Q,P_{2n},\tau)$ is a modular form on $\Gamma_0(N)$ of weight $r+2n$ 
and character $\e.$ This is the theorem of Schoeneberg [S1].

Suppose next that $v$ is a unit vector. Then the transformation law
(\ref{1.4}) says that $\Theta(Q,v,\tau,X)$ is a (holomorphic) 
Jacobi-like form of weight $r,$ level $N,$ character $\e$ in the
sense of Zagier [Z]. 

When $v$ is a unit vector, $\theta(Q,v,k,\t)$ will not be  modular,
but one can suitably combine two Jacobi-like forms to produce a
sequence of modular forms (loc. cit.). For example, there is a
holomorphic Jacobi-like form of weight $0$  of particular interest,
namely
\begin{equation}\label{1.6}
\tilde E_2(\t,X)=\sum_{n\geq 0}(-1)^n\frac{E_2(\tau)^n}{n!}(2\pi iX)^n
\end{equation}
where 
$$E_2(\t)=-\frac{1}{12}+2\sum_{n\geq 1}\sigma_1(n)q^n$$
is the
usual ``unmodular'' Eisenstein series of ``weight'' 2 normalized as
indicated. Using this together with Theorem \ref{t1} yields
 \begin{th}\label{t2} Let notation be as before, and suppose that $v$ is 
a unit vector. Set
\begin{equation}\label{1.7}
\gamma(t,k)=2^{-t}{k\choose t}{k-t\choose t}t!
\end{equation}
\begin{equation}\label{1.8}
\Psi(Q,v,2k,\tau)=\sum_{t=0}^{k}\gamma(t,2k)E_2(\tau)^t\theta(Q,v,2k-2t,\t).
\end{equation}
Then $\Psi(Q,v,2k,\tau)$ is a holomorphic modular form on $\Gamma_0(N)$ 
 of weight $2k+r$ and character $\e.$
\end{th}

One knows that if $k>0$ and $v$ is a null vector then in fact
$\theta(Q,P_k,\tau)$ is a cusp form. In the same spirit we have the 
following supplement to Theorem \ref{t2}:
\begin{th}\label{t3} Let notation and assumptions be as in Theorem
\ref{t2}, and $k\geq 2.$ Then 
\begin{equation}\label{1.9}
\Psi(Q,v,2k,\tau)-\gamma(k,2k)(-\frac{1}{12})^{k} \theta(Q,\tau)E_{2k}(\tau)
\end{equation}
is a cusp form on $\Gamma_0(N)$ of weight $2k+r$ and character
$\e.$ Here, 
$$E_{2k}(\tau)=1-\frac{4k}{B_{2k}}\sum_{n\geq 1}\sigma_{k-1}(n)q^n$$
is the usual Eisenstein series.
\end{th}

It is forms of the general shape (\ref{1.9}) that appear as partition function
in certain conformal field theories [DMN], and whose existence led us to 
the results of the present paper.

The paper is organized as follows: in Section 2 we discuss Jacobi-like
forms and the proof of Theorem \ref{t2}. Section 3 is devoted to the
proof of Theorems \ref{t1} and \ref{t3}, which follows in general
outline the original proof of Schoeneberg. In Section 4 we discuss
connections with conformal field theory and state some further results which
will be proved in [DMN].

\section{Jacobi-like forms}
\setcounter{equation}{0}
We are interested in holomorphic functions $\phi(\t,X)$ on ${\frak
h}\times \C$ of the form 
\begin{equation}\label{2.1}
\phi(\t,X)=\sum_{n\geq 0}\phi^{(n)}(\tau)(2\pi i X)^n
\end{equation}
and which satisfy 
\begin{equation}\label{2.2}
\phi(\gamma\t,\frac{X}{(c\tau+d)^2})=\chi(d)(c\tau+d)^k\exp\left(\frac{cmX}{c\tau+d}\right)\phi(\tau,X)
\end{equation}
for some $m\in \C,$ integer $k$ and Dirichlet character $\chi (\mod N)$, 
and for all $\g\in \left(\begin{array}{cc} a & b\\ c & d\end{array}\right)\in
\Gamma_0(N).$ By scaling  $X,$ the essential cases correspond to $m=0$ and 
$m=1.$ As long as $\phi$ is holomorphic
at the cusps, the case $m=0$ means precisely 
that each $\phi^{(n)}(\tau)$ is a holomorphic modular form
on $\Gamma_0(N)$ of weight $k+2n$ and character $\chi.$ The case $m=1$ means
that $\phi$ is a holomorphic Jacobi-like form on $\Gamma_0(N)$ of weight 
$k$ and character $\chi$ (cf. [Z]).

By way of examples, let $\tilde E_2(\t,X)$ be as in (\ref{1.6}). The particular normalization
of $E_2(\t)$ that we are using satisfies the well-known transformation
law
\begin{equation}\label{2.3}
E_2(\gamma\tau)=(c\tau+d)^2E_2(\tau)-\frac{c(c\tau+d)}{2\pi i}.
\end{equation}

It follows that $\tilde E_2(\tau,X)=\exp(-2\pi i E_2(\t)X)$
satisfies 
\begin{eqnarray*}
& &\tilde E_2(\gamma\t,\frac{X}{(c\tau+d)^2})=\exp\left(-2\pi i\left((c\tau+d)^2E_2(\tau)-\frac{c(c\tau+d)}{2\pi i}\right)\frac{X}{(c\t+d)^2}\right)\\
& &\ \ \ \ \ \ \ \ \ \ \ \ \ =\exp\left(\frac{cX}{c\t+d}\right)\tilde E_2(\tau,X),
\end{eqnarray*}
so that $\tilde E_2(\tau,X)$ is a holomorphic Jacobi-like form
of level 1 and weight 0.

Now let $\Theta(\!Q,v,\tau,X\!)$ be as in (\ref{1.3}) with $v$ a unit
vector $v.$ Then $\tilde E_2(\tau,-X)\Theta(Q,v,\tau,X)$ satisfies
(\ref{2.2}) with $m=0$ (we are assuming the truth of Theorem \ref{t1} at this point). We have
$$\tilde E_2(\tau,-X)\Theta(Q,v,\tau,X)=\sum_{k\geq
0}f^{(k)}(\tau)(2\pi iX)^k$$ 
where 
\begin{equation}\label{2.4}
f^{(k)}(\tau)=\frac{2^k}{(2k)!}\sum_{t=0}^{k}\gamma(t,2k)E_2(\tau)^t\theta(Q,v,2k-2t,\t)
\end{equation}
and $\gamma(t,2k)$ is as in (\ref{1.7}). It follows, granting
holomorphy at the cusps for now, that $f^{(k)}(\tau)$ is a holomorphic
modular form on $\Gamma_0(N)$ of weight $r+2k$ and character $\e.$ So Theorem \ref{t2} follows from Theorem \ref{t1}.

There are other Jacobi-like forms that one could use in place of
$\tilde E_2(\tau,X)$ in order to construct modular forms involving the
$\theta(Q,v,n,\tau).$ For example, we could take the Cohen-Kuznetsov
Jacobi-like form 
$$\sum_{n\geq 0}\frac{f^{(n)}(\tau)}{n!(n+k-1)!}(2\pi i X)^n$$
for $f$ a modular form of weight $k,$ as described in [Z]. We will not pursue this possibility
here: it is the forms $\Psi(Q,v,2k,\t)$ of Theorem \ref{t2} that we need in 
[DMN].

\section{Proofs of Theorem 1 and 3}
\setcounter{equation}{0}
In this section we present the proofs of Theorem \ref{t1} and
\ref{t3}. They follow in general outline the proof of Schoeneberg
[S2]. We therefore adopt notation similar to (loc. cit.) and omit some
details. In particular, we have $\<x,y\>=x'Ay$ for $x,y\in \C^f,$
where $x'$ denotes the transpose of the column vector $x.$

Let $A$ and $Q$ be as in Section 1. For $x=(x_1,...,x_f)$ and a scalar
$\l$ we set
\begin{equation}\label{3.1}
\theta_{\l}(A,x)=\sum_{m\in\Z^f}\exp(2\l Q(m+x)).
\end{equation}

For $l=(l_1,...,l_f)\in \C^f$ we let ${\cal L}$ be the linear differential
operator
\begin{equation}\label{3.2}
{\cal L}=\sum_{i=1}^fl_i\frac{\partial}{\partial x_i}.
\end{equation}

\begin{lem}\label{l3.1} Let $k\geq 0$ be an integer. Then 
\begin{equation}\label{3.3}
{\cal L}^k(\theta_{\l}(A,x))=\sum_{i=0}^{[k/2]}\sum_{m\in\Z^f}\gamma(i,k)(2\l)^{k-i}(2Q(l))^i(l'A(m+x))^{k-2i}\exp(2\l Q(m+x))
\end{equation}
where $\gamma(i,k)$ is defined by (\ref{1.7}).
\end{lem}

\pf One sees that there is an equality of the form
\begin{equation}\label{3.4}
{\cal L}^k(\theta_{\l}(A,x))=\sum_{i=0}^{[k/2]}\sum_{m\in\Z^f}\gamma_{\l}(i,k)
(2Q(l))^i(l'A(m+x))^{k-2i}\exp(2\l
Q(m+x))
\end{equation}
for some scalars $\gamma_{\l}(i,k), 0\leq i\leq [k/2], k\geq 0.$
Setting $\gamma_{\l}(i,k)=0$ for values of $i$ and $k$ not in these
ranges, $\gamma_{\l}(i,k)$ satisfies a recursion relation, namely
\begin{equation}\label{3.5}
\gamma_{\l}(i,k+1)=(k+2-2i)\gamma_{\l}(i-1,k)+2\l\gamma_{\l}(i,k),
\gamma_{\l}(0,0)=1.
\end{equation}

We can solve the recursion, and find that $\gamma_{\l}(i,k)=(2\l)^{k-i}\gamma(i,k).$ The lemma follows. \qed

Now one knows (e.g. page 205 of [S2]) that the following
transformation law holds: for $\t$ in the upper half-plane and
for a suitable determination of the square root,
\begin{eqnarray}\label{3.6}
& &\ \ \ \ \sum_{m\in\Z^f}\exp(2\pi i \t Q(m+x))=\theta_{\pi i\tau}(A,x)\nonumber\\
& & =
\frac{1}{(\sqrt{-i\t})^f(\det A)^{1/2}}\sum_{m\in\Z^f}\exp(-\frac{\pi i}{\t}m'A^{-1}m+2\pi i m'x).
\end{eqnarray}

We apply the operator ${\cal L}^k$ to both sides of (\ref{3.6}), using Lemma
\ref{l3.1}, to obtain 
\begin{eqnarray}\label{3.7}
& &\ \ \ \ \sum_{j=0}^{[k/2]}\sum_{m\in\Z^f}(2\pi i \t)^{k-j}\gamma(j,k)(2Q(l))^j
(l'A(m+x))^{k-2j}\exp(2\pi i \t Q(m+x))\nonumber\\
& &=\frac{1}{(\sqrt{-i\t})^f(\det A)^{1/2}}\sum_{m\in\Z^f}(2\pi im'l)^k
\exp(-\frac{\pi i}{\t}m'A^{-1}m+2\pi i m'x).
\end{eqnarray}

Recall that $N$ is the level of $A.$ Following [S2], we replace $x$ by $h/N,$
$\t$ by $-1/\t$ and $m$ by $Am_1/N$  on the r.h.s. of (\ref{3.7}). Remembering that $f=2r$ and setting 
$D=\det A,$ we get 
\begin{eqnarray}\label{3.8}
& &\ \ \ \ \frac{(2\pi i )^k\t^rN^{-k}}{i^r\sqrt{D}}\sum_{\stackrel{m_1\in\Z^f}{Am_1=0 (N)}}(l'Am_1)^{k}\exp(2\pi i \t Q(m_1)/N^2+2\pi i
m_1Ah/N^2)\nonumber\\ 
& &=\sum_{j=0}^{[k/2]}\sum_{\stackrel{m\in\Z^f}{ m\equiv h (N)}}N^{2j}\left(\frac{-2\pi i}{\t}\right)^{k-j}\gamma(j,k)(2Q(l))^j
(l'Am)^{k-2j}\exp(\frac{-2\pi i}{\t}\frac{Q(m)}{N^2}).
\end{eqnarray}

Note that if $l$ is a null vector, only the term with $j=0$ survives 
on the r.h.s. of (\ref{3.8}), which then reduces to equation
(12) of [S2], page 209.

\bigskip
 
We discuss some transformation formulas. 
For $h\in\Z^f,$ $Ah\equiv 0\ (\mod N)$ and $0\leq j\leq k,$ we set 
\begin{equation}\label{3.9}
\theta(A,h,l,k,\t)=\frac{1}{N^k}\sum_{\stackrel{m\in \Z^f}{m\equiv h (N)}}
(l'Am)^k\exp(2\pi i\tau Q(m)/N^2).
\end{equation}
It is also convenient to introduce
\begin{equation}\label{3.10} 
\Theta(A,h,l,k,j,\t)=\frac{(-i)^{r+2k}\t^{r+k-2j}}{\sqrt{D}}
\sum_{\stackrel{g\,\mod N}{Ag\equiv 0 (N)}}\exp(2\pi i g'Ah/N^2)\theta(A,g,l,k-2j,\tau).
\end{equation}

\begin{rem}\label{rem}{\em  (i) It is clear that $\Theta(A,h,l,k,j,\t)$ depends only on 
$k-2j,$ rather than both $k$ and $j.$ However it is convenient to keep 
the notation as it is.

(ii) Note that if $h=0$ then $\theta(A,h,l,k,\t)$ is just the function $\theta(A,l,k,\tau)=\theta(Q,l,k,\tau).$}
\end{rem}

\begin{thm}\label{t3.2} We have 
\begin{equation}\label{3.11} 
\theta(A,h,l,k,-1/\t)=\sum_{j=0}^{[k/2]}\left(\frac{Q(l)\tau}{\pi i}\right)^j
\gamma(j,k)\Theta(A,h,l,k,j,\t).
\end{equation}
\end{thm}

\pf As $Ah\equiv 0\ \mod N,$ we can split-off an exponential factor form
the l.h.s. of (\ref{3.8}). Then using (\ref{3.9}), (\ref{3.8}) can be
written in the form 
\begin{eqnarray}\label{3.12}
& &\ \ \ \ \frac{(2\pi i )^k\t^r}{i^r\sqrt{D}}\sum_{\stackrel{g\,\mod N}{ Ag\equiv 0 (N)}}
\exp(2\pi i g'Ah/N^2)\theta(A,g,l,k,\tau)\nonumber\\
& &=\sum_{j=0}^{[k/2]}\left(\frac{-2\pi i}{\t}\right)^{k-j}
\gamma(j,k)(2Q(l))^j\theta(A,h,l,k-2j,-1/\t).
\end{eqnarray}

We shall prove Theorem \ref{t3.2} by induction on $k.$ If $k=0$ it reduces to a standard transformation law (equation 17II of [S2], page 210). In the general
case, (\ref{3.12}) yields
\begin{eqnarray*}
& &\theta(A,h,l,k,-1/\t)=\frac{(-i )^{r+2k}\t^{r+k}}{\sqrt{D}}\sum_{\stackrel{g\,\mod N}{ Ag\equiv 0 (N)}}\exp(2\pi i g'Ah/N^2)\theta(A,g,l,k,\tau)\\
& &\ \ \ \ \ \ \ \ \ \ \ \ \ \ \ \  \ \ -
\sum_{j=1}^{[k/2]}\gamma(j,k)\left(-\frac{Q(l)\tau}{\pi i}\right)^j
\theta(A,h,l,k-2j,-1/\t)\\
& &=\Theta(A,h,l,k,0,\t)\\
& & \ \ \ \ -\sum_{j=1}^{[k/2]}\gamma(j,k)\left(\frac{-Q(l)\t}{\pi i}\right)^j\sum_{t=0}^{[(k-2j)/2]}\gamma(t,k-2j)\left(\frac{Q(l)\tau}{\pi i}\right)^t\Theta(A,h,l,k-2j,t,\t).
\end{eqnarray*}

As $\Theta(A,h,l,k-2j,t,\t)=\Theta(A,h,l,k,t+j,\t),$ we see that 
$\theta(A,h,l,k,-1/\t)$ is equal to 
\begin{equation}\label{3.13}
\Theta(A,h,l,k,0,\t)+\sum_{u=1}^{[k/2]}\left(\frac{Q(l)\tau}{\pi i}\right)^u
\beta(u,k)\Theta(A,h,l,k,u,\t)
\end{equation}
where 
\begin{equation}\label{3.14}
\beta(u,k)=-\sum_{j=1}^{[k/2]}(-1)^j\gamma(j,k)\gamma(u-j,k-2j).
\end{equation}
From (\ref{1.7}) and (\ref{3.14}) we see that 
$$\beta(u,k)=\gamma(u,k)-\sum_{j=0}^{[k/2]}(-1)^j\gamma(j,k)\gamma(u-j,k-2j)
=\gamma(u,k)-\sum_{j=0}^{[k/2]}(-1)^j\gamma(u,k){u\choose j}$$
i.e., $\beta(u,k)=\gamma(u,k).$ Now (\ref{3.13}) implies the desired equality
(\ref{3.11}). \qed

\bigskip

We now proceed to our main transformation formula, which is the following:
\begin{thm}\label{t3.3} If
$\left(\begin{array}{cc} a & b\\ c &d\end{array}\right)$ lies in $\Gamma_0(N),$ and
if $d>0,$ then 
\begin{eqnarray}\label{3.15}
& &\ \ \ \ (c\t+d)^{-(r+k)}\theta(A,h,l,k,\frac{a\t+b}{c\t+d})\nonumber\\
& &=\exp(2\pi i Q(h)ab/N^2)\e(d)\sum_{j=0}^{[k/2]}\left(\frac{Q(l)c}{\pi i(c\tau +d)}\right)^j\gamma(j,k)\theta(A,bh,l,k-2j,\t).
\end{eqnarray}
In particular, taking $h=0,$ if $d>0$ then 
\begin{equation}\label{3.16}
(c\t+d)^{-(r+k)}\theta(A,l,k,\frac{a\t+b}{c\t+d})
=\e(d)\sum_{j=0}^{[k/2]}\left(\frac{Q(l)c}{\pi i(c\tau +d)}\right)^j\gamma(j,k)\theta(A,l,k-2j,\t).
\end{equation}
\end{thm}

We begin by noting that
\begin{equation}\label{3.17}
\theta(A,h,l,k,\t+1)=\exp(2\pi iQ(h)/N^2)\theta(A,h,l,k,\tau),
\end{equation}
and also if $c>0$ then 
\begin{equation}\label{3.18}
\theta(A,h,l,k,\t)=\sum_{\stackrel{g\,\mod cN}{ g\equiv h (N)}}\theta(cA,g,l,k,c\tau).
\end{equation}
(\ref{3.17}) is immediate from (\ref{3.9}) (remembering that $Ah=0\ (\mod N)$);
(\ref{3.18}) follows as in equation 18 of [S2], page 211.

Using (\ref{3.17}), (\ref{3.18}) and Theorem \ref{t3.2} we calculate
for $\gamma=\left(\begin{array}{cc} a & b\\ c &d\end{array}\right)$
in $SL(2,\Z)$ with $c>0$ that 
\begin{eqnarray}\label{3.19}
& & \theta(A,h,l,k,\gamma\t)=\theta(A,h,l,k,c^{-1}(a-(c\tau+d)^{-1}))\nonumber\\
& &=\sum_{\stackrel{g\,\mod cN}{g\equiv h (N)}}\theta(cA,g,l,k,a-(c\tau+d)^{-1})\nonumber\\
& &=\sum_{\stackrel{g\,\mod cN}{ g\equiv h (N)}}\exp(2\pi iacQ(g)/c^2N^2)\theta(cA,g,l,k,
-(c\tau+d)^{-1})\nonumber\\
& &=\sum_{\stackrel{g\,\mod cN}{ g\equiv h (N)}}\sum_{j=0}^{[k/2]}\exp(2\pi iaQ(g)/cN^2)
\left(\frac{Q(l)c(c\t+d)}{\pi i}\right)^j\gamma(j,k)\Theta(cA,g,l,k,j,c\tau+d)\nonumber\\
& &=\sum_{\stackrel{g\,\mod cN}{ g\equiv h (N)}}\sum_{j=0}^{[k/2]}\sum_{\stackrel{
q\,\mod cN}{cAq\equiv 0 (cN)}}\exp(2\pi iaQ(g)/cN^2)\left(\frac{Q(l)c(c\t+d)}{\pi i}\right)^j
\gamma(j,k)\frac{(-i)^{r+2(k-2j)}}{\sqrt{c^fD}}\cdot\nonumber\\
& &\ \ \ \ \ \ \ \ \ \ \ \ \ \ \ \ \cdot(c\tau+d)^{r+k-2j}\exp(2\pi ig'Aq/cN^2)\theta(cA,q,l,k-2j,c\tau+d)\nonumber\\
& &=\frac{(c\tau+d)^{r+k}(-i)^{r+2k}}{c^r\sqrt{D}}\sum_{\stackrel{g\,\mod cN}{ g\equiv h (N)}}\sum_{\stackrel{q\,\mod cN}{ Aq\equiv 0 (N)}}\sum_{j=0}^{[k/2]}\exp(2\pi i(aQ(g)+dQ(q)+g'Aq)/cN^2)\cdot\nonumber\\
& &\ \ \ \ \ \ \ \ \ \ \ \ \ \ \ \ \ \cdot\gamma(j,k)\left(\frac{cQ(l)}{(c\t+d)\pi i}\right)^j\theta(cA,q,l,k-2j,c\t)
\end{eqnarray}
where we have used the fact that $cA$ has level $cN$.

Following [S2], page 213 we write
\begin{equation}\label{3.20}
\phi_{h,q}=\sum_{\stackrel{g\,\mod cN}{ g\equiv h (N)}}\exp(2\pi i(aQ(g)+dQ(q)+g'Aq)/cN^2)
\end{equation}
and note that $\phi_{h,q}$ depends on $q$ only modulo $N.$ Then 
(\ref{3.19}) can be put into the form
\begin{eqnarray*}
& & \ \ \ \ \frac{(c\tau+d)^{r+k}(-i)^{r+2k}}{c^r\sqrt{D}}\sum_{j=0}^{[k/2]}
\sum_{\stackrel{q\,\mod cN}{ Aq\equiv 0 (N)}}\phi_{h,q}\gamma(j,k)
\left(\frac{cQ(l)}{(c\t+d)\pi i}\right)^j\theta(cA,q,l,k-2j,c\t)\\
& &=\frac{(c\tau+d)^{r+k}(-i)^{r+2k}}{c^r\sqrt{D}}\sum_{j=0}^{[k/2]}
\sum_{q_1\,n \mod N} \phi_{h,q_1}\gamma(j,k)\left(\frac{cQ(l)}{(c\t+d)\pi i}\right)^j
\cdot\\
& &\ \ \ \ \ \ \ \ \ \ \ \ \ \ \ \ \cdot\sum_{\stackrel{q\,\mod cN}{Aq\equiv0 (N), q\equiv q_1 (N)}}\theta(cA,q,l,k-2j,c\t).
\end{eqnarray*}

Another application of (\ref{3.18}) now yields 
\begin{lem}\label{l3.4} Let 
$\gamma=\left(\begin{array}{cc} a & b\\ c &d\end{array}\right)\in SL(2,\Z)$ with $c>0.$ Then 
\begin{eqnarray}\label{3.21}
& &\ \ \ \ (c\tau+d)^{-(r+k)}\theta(A,h,l,k,\gamma\tau)\nonumber\\
& &=\frac{(-i)^{r+2k}}{c^r\sqrt{D}}\sum_{j=0}^{[k/2]}
\sum_{\stackrel{q_1\,\mod N}{ Aq_1\equiv 0 (N)}} \phi_{h,q_1}\gamma(j,k)
\left(\frac{cQ(l)}{(c\t+d)\pi i}\right)^j\theta(A,q_1,l,k-2j,\t).
\end{eqnarray}
\end{lem}

We continue the calculation, but now assuming also that $d\equiv 0 \ (\mod N).$
As in (loc. cit.) we now find that (\ref{3.21}) can be written as follows:
\begin{eqnarray}\label{3.22}
& &(c\tau+d)^{-(r+k)}\theta(A,h,l,k,\gamma\tau)=
\frac{(-i)^{r+2k}\phi_{h,0}}{c^r\sqrt{D}}\sum_{j=0}^{[k/2]}
\sum_{\stackrel{q_1\,\mod N}{Aq_1\equiv 0 (N)}}\exp(-2\pi i h'Aq_1b/N^2)\cdot\nonumber\\
& &\ \ \ \ \ \ \ \ \ \ \ \ \ \ \ \cdot\gamma(j,k)
\left(\frac{cQ(l)}{(c\t+d)\pi i}\right)^j\theta(A,q_1,l,k-2j,\t).
\end{eqnarray}

Replacing $\t$ by $-1/\t$ in (\ref{3.22}) and using Theorem
\ref{t3.2} again leads to
\begin{eqnarray*}
& &\ \ \ \ \left(\frac{d\tau-c}{\tau}\right)^{-(r+k)}\theta(A,h,l,k,\frac{b\tau-a}{d\tau-c})\\
& &=\frac{(-i)^{r+2k}\phi_{h,0}}{c^r\sqrt{D}}\sum_{j=0}^{[k/2]}
\sum_{\stackrel{q_1\,\mod N}{ Aq_1\equiv 0 (N)}}\sum_{u=0}^{[(k-2j)/2]}\exp(-2\pi i h'Aq_1b/N^2)\gamma(j,k)\cdot\\
& &\ \ \ \ \ \ \ \ \ \ \ \ \cdot\left(\frac{cQ(l)\t}{(d\t-c)\pi i}\right)^j\left(\frac{Q(l)\tau}{\pi i}\right)^u\gamma(u,k-2j)\Theta(A,q_1,l,k-2j,u,\t)\\
& &=\frac{(-i)^{r+2k}\phi_{h,0}}{c^r\sqrt{D}}\sum_{j=0}^{[k/2]}
\sum_{\stackrel{q_1\,\mod N}{ Aq_1\equiv 0 (N)}}\sum_{u=0}^{[(k-2j)/2]}
\sum_{\stackrel{g\,\mod N}{Ag\equiv 0 (N)}}\exp(-2\pi i h'Aq_1b/N^2)\cdot\\
& &\ \ \ \ \ \ \ \ \ \ \ \ \ \cdot\gamma(j+u,k)
{j+u\choose j}\left(\frac{c}{d\t-c}\right)^j
\left(\frac{Q(l)\tau}{\pi i}\right)^{j+u}
\frac{(-i)^{r+2(k-2j)}}{\sqrt{D}}\cdot\\
& &\ \ \ \ \ \ \ \ \ \ \ \ \ \cdot\t^{r+k-2j-2u}\exp(2\pi i g'Aq_1/N^2)
\theta(A,g,l,k-2j-2u,\t).
\end{eqnarray*}

Therefore, we get 
\begin{eqnarray}\label{3.23}
& &\ \ \ \ (d\tau-c)^{-(r+k)}\theta(A,h,l,k,\frac{b\tau-a}{d\tau-c})\nonumber\\
& &=\frac{(-1)^{r}\phi_{h,0}}{c^rD}\sum_{j=0}^{[k/2]}\sum_{\stackrel{q_1\,\mod N}{ Aq_1\equiv 0 (N)}}\sum_{u=0}^{[(k-2j)/2]}
\sum_{\stackrel{g\,\mod N}{ Ag\equiv 0 (N)}}\exp(2\pi i (g-bh)'Aq_1/N^2)
\gamma(j+u,k)\cdot\nonumber\\
& &\ \ \ \ \ \ \ \ \  \ \ \ \cdot{j+u\choose j}\left(\frac{c}{d\t-c}\right)^j
\left(\frac{Q(l)}{\pi i\tau}\right)^{j+u}
\theta(A,g,l,k-2(j+u),\tau).
\end{eqnarray}

Now we know  (page 214 of [S2]) that 
$$\sum_{\stackrel{q_1\,\mod N}{Aq_1\equiv 0 (N)}}\exp(2\pi i(g-bh)'Aq_1/N^2)=D\delta_{g,bh}$$
where $\delta_{g,bh}$ is the Kronecker delta and $g,bh$ are considered modulo
$N.$ Thus (\ref{3.23}) reduces to
 \begin{eqnarray*}
& &(d\tau-c)^{-(r+k)}\theta(A,h,l,k,\frac{b\tau-a}{d\tau-c})\nonumber\\
& &=\frac{(-1)^{r}\phi_{h,0}}{c^r}
\sum_{j=0}^{[k/2]}\sum_{u=0}^{[(k-2j)/2]}\gamma(j+u,k){j+u\choose j}\left(\frac{c}{d\t-c}\right)^j
\left(\frac{Q(l)}{\pi i\tau}\right)^{j+u}\cdot\\
& &\ \ \ \ \ \ \ \ \ \ \  \ \ \ \ \cdot\theta(A,bh,l,k-2(j+u),\tau)\\
& &=\frac{(-1)^{r}\phi_{h,0}}{c^r}
\sum_{j=0}^{[k/2]}\sum_{t=0}^{[k/2]}\gamma(t,k){t\choose j}
\left(\frac{c}{d\t-c}\right)^j
\left(\frac{Q(l)}{\pi i\tau}\right)^{t}
\theta(A,bh,l,k-2t,\tau)\\
& &=\frac{(-1)^{r}\phi_{h,0}}{c^r}\sum_{t=0}^{[k/2]}\gamma(t,k)\left(\frac{Q(l)}{\pi i\tau}\right)^{t}\left(\frac{d\t}{d\t-c}\right)^t\theta(A,bh,l,k-2t,\tau).
\end{eqnarray*}

Making the change of variables 
$\left(\begin{array}{cc} b & -a\\ d &-c\end{array}\right)\mapsto 
\left(\begin{array}{cc} a & b\\ c &d\end{array}\right)$ yields the equality
\begin{eqnarray}\label{3.24}
& &\ \ \  \ (c\tau+d)^{-(r+k)}\theta(A,h,l,k,\frac{a\tau+b}{c\t+d})\nonumber\\
& &=\frac{\phi_{h,0}}{d^r}\sum_{t=0}^{[k/2]}\gamma(t,k)\left(\frac{Q(l)c}{\pi i(c\tau+d)}\right)^{t}\theta(A,bh,l,k-2t,\tau).
\end{eqnarray}

Finally, it is known (page 215 et seq of [S2]) that
$$ \frac{\phi_{h,o}}{d^r}=\exp(2\pi i Q(h)ab/N^2)\e(d).$$
Now (\ref{3.24}) implies all assertions of Theorem \ref{t3.3}.

\bigskip

We can now complete the proof of Theorem \ref{t1}. Noting Remark \ref{rem} (ii)
and using (\ref{3.16}) with
$\gamma=\left(\begin{array}{cc} a & b\\ c &d\end{array}\right),$ we get
\begin{eqnarray*} 
& &\Theta(Q,v,\gamma\tau,\frac{X}{(c\tau+d)^2})=\sum_{n\geq 0}
\frac{2^n\theta(Q,v,2n,\gamma\tau)}{(2n)!}\left(\frac{2\pi i X}{ (c\tau+d)^2}
\right)^n\\
& &=\e(d)(c\tau+d)^r\sum_{n\geq 0}\sum_{j=0}^n
\frac{2^n}{(2n)!}\left(\frac{\<v,v\>c}{2\pi i(c\tau+d)}\right)^j\gamma(j,2n)
(c\tau +d)^{2n}\cdot\\
& &\ \ \ \ \ \ \ \ \ \ \ \cdot \theta(Q,v,2n-2j,\tau)\left(\frac{2\pi i X}{ (c\tau+d)^2}
\right)^n\\
& &=\e(d)(c\tau+d)^r\sum_{n\geq 0}\sum_{j=0}^n
\frac{2^{n-j}}{(2n-2j)!}\theta(Q,v,2n-2j,\tau)(2\pi iX)^{n-j}
\left(\frac{\<v,v\>c}{c\tau+d}\right)^j\frac{X^j}{j!}\\
& &=\e(d)(c\tau+d)^r\exp\left(\frac{\<v,v\>cX}{c\tau+d}\right)
\Theta(Q,v,\tau,X),
\end{eqnarray*}
which is the required (\ref{1.4}) in the case $d>0.$ The general
case follows easily.

\bigskip

Finally, we consider behavior at the cusps. This will lead to the proof of Theorem \ref{t3} as well as the proof of Theorem \ref{t2} initiated in Section 2.

To check the expansion of $\Psi(Q,l,k,\tau)$ at the finite cups, we use Lemma \ref{l3.4}. Thus for $c>0,$
\begin{eqnarray}\label{3.25}
& &\ \ \ \ (c\tau+d)^{-(r+k)}\Psi(Q,l,k,\frac{a\tau+b}{c\t+d})\nonumber\\
& &=\sum_{t=0}^{[k/2]}\gamma(t,k)(c\tau+d)^{-2t}E_2(\frac{a\t+b}{c\t+d})^t
(c\tau+d)^{-(r+k-2t)}\theta(A,l,k-2t,\frac{a\tau+b}{c\t+d})\nonumber\\
& &=\sum_{t=0}^{[k/2]}\gamma(t,k)\left(E_2(\t)-\frac{c}{2\pi i(c\tau+d)}\right)^t\frac{(-i)^{r+2(k-2t)}}{c^r\sqrt{D}}
\sum_{j=0}^{[(k-2t)/2]}\sum_{\stackrel{q\,\mod N}{Aq\equiv 0 (N)}}\phi_{0,q}\cdot\nonumber\\
& &\ \ \ \ \ \ \ \ \ \  \cdot\gamma(j,k-2t)\left(\frac{cQ(l)}{\pi i(c\tau+d)}\right)^j
\theta(A,q,l,k-2t-2j,\tau)\nonumber\\
& &=\frac{(-i)^{r+2k}}{c^r\sqrt{D}}\sum_{t=0}^{[k/2]}\sum_{j=0}^{[(k-2t)/2]}
\sum_{\stackrel{q\,\mod N}{ Aq\equiv 0 (N)}}\phi_{0,q}\gamma(j+t,k)
{j+t\choose j}\left(\frac{c}{2\pi i(c\tau+d)}\right)^j\cdot\nonumber\\
& &\ \ \ \ \ \ \ \ \ \ \ \cdot(2Q(l))^j
\left(E_2(\t)-\frac{c}{2\pi i(c\tau+d)}\right)^t\theta(A,q,l,k-2t-2j,\tau)\nonumber\\
& &=
\frac{(-i)^{r+2k}}{c^r\sqrt{D}}\sum_{\stackrel{q\,\mod N}{ Aq\equiv 0 (N)}}\phi_{0,q}
\sum_{s=0}^{[k/2]}\gamma(s,k)E_2(\tau)^s\theta(A,q,l,k-2s,\tau).
\end{eqnarray}
This shows that $\Psi(Q,l,k,\tau)$ is holomorphic at the cusps, and
thus complete the proof of Theorem \ref{t2}.

As for Theorem \ref{t3}, the value of $\Psi(Q,v,2k,\tau)$ at $i\infty$
is the same as that of the function
$$\gamma(k,2k)E_2(\tau)^k\theta(A,v,0,\tau)=\gamma(k,2k)E_2(\tau)^k\theta(Q,\tau)$$
which has $q$-expansion $\gamma(k,2k)(-1/12)^k(1+\cdots)$. Thus
$$\Psi(Q,v,2k,\tau)-\gamma(k,2k)(-1/12)^k\theta(Q,\tau)E_{2k}(\tau)$$
certainly vanishes at $i\infty.$ From (\ref{3.25}), the value at a general cusp is
$$\frac{(-i)^{r}}{c^r\sqrt{D}}\sum_{\stackrel{q\,\mod N}{ Aq\equiv 0 (N)}}\phi_{0,q}\gamma(k,2k)(-1/12)^k.$$

But $\displaystyle{\frac{(-i)^{r}}{c^r\sqrt{D}}\sum_{\stackrel{q\,\mod N}{ Aq\equiv 0 (N)}}\phi_{0,q}}$ is the value of $\theta(Q,\tau)$ at the same cusp 
(cf. equation (24) of [S2], page 213), whence it is clear that $\Psi(Q,v,2k,\tau)
-\gamma(k,2k)(-1/12)^k\theta(Q,\tau)E_{2k}(\tau)$ vanishes at every cusp if $k\geq 2.$ This completes the proof of Theorem \ref{t3}.

\section{Concluding comments}
\setcounter{equation}{0}
We have already mentioned that the previous results were motivated by conformal 
field theory, more precisely by the problem of
calculating $1$-point correlation functions for vertex operator algebras
[DMN]. For earlier results in this direction, see [DLM] and [DM]. This
perspective also enables us to prove the following result (see [DMN]):
let the notation be as before, and suppose that the 
lattice $\Z^f$ contains a root $\alpha$ ie., $Q(\alpha)=1.$ Then the cusp form of Theorem \ref{t3} (with $v=\frac{\alpha}{\sqrt{2}}$) is identically zero. 
That is, we have
\begin{equation}\label{4.1}
\theta(Q,\frac{\alpha}{\sqrt{2}},4,\tau)+6E_2(\tau)\theta(Q,\frac{\alpha}{\sqrt{2}},2,\tau)+3E_2(\tau)^2\theta(Q,\tau)=\frac{1}{48}E_4(\tau)\theta(Q,\tau).
\end{equation}

It is interesting that in [Z], Zagier raises the question of whether
there is a relation between Jacobi-like forms and vertex operator
algebras. The present work together with [DMN] certainly suggests that
this question continues to be one which is worth exploring.

\end{document}